\documentclass[reqno, table]{amsart} 
\usepackage{amsmath}
\usepackage[margin=1.0in]{geometry}
\usepackage{amscd,amsmath,amsxtra,amsthm,amssymb,stmaryrd,xr,mathrsfs,mathtools, enumerate,commath, comment, enumitem, graphicx}
\usepackage{amsfonts, amssymb, amsthm, amsmath, braket, hyperref, mathtools, comment, mathrsfs, enumitem, cleveref}
\usepackage{stmaryrd}
\usepackage{amsfonts}
\usepackage{orcidlink}
\usepackage{multirow}
\usepackage{xcolor}
\usepackage{commath}
\usepackage{comment}
\usepackage{tikz-cd}
\usepackage{tkz-graph}
\usepackage{colonequals}
\usepackage{hyperref}
\usepackage{graphicx}
\usepackage{bigstrut}
\usepackage{enumerate}
\usepackage{tabularx}
\numberwithin{equation}{section}
\usepackage{longtable, multirow, array}
\usepackage{cleveref}
\crefname{section}{§}{§§}
\Crefname{section}{§}{§§}
\newcommand{\genlegendre}[4]{%
	\genfrac{(}{)}{}{#1}{#3}{#4}%
	\if\relax\detokenize{#2}\relax\else_{\!#2}\fi
}
\newcommand{\legendre}[3][]{\genlegendre{}{#1}{#2}{#3}}

\usepackage{longtable} 
\usepackage{pdflscape} 
\usepackage{booktabs}
\usepackage{hyperref}
\definecolor{vegasgold}{rgb}{0.77, 0.7, 0.35}
\definecolor{darkgoldenrod}{rgb}{0.72, 0.53, 0.04}
\definecolor{gold(metallic)}{rgb}{0.83, 0.69, 0.22}
\hypersetup{
 colorlinks=true,
 linkcolor=darkgoldenrod,
 filecolor=brown,      
 urlcolor=gold(metallic),
 citecolor=darkgoldenrod,
 }
\newtheorem{lthm}{Theorem}

\usepackage[all,cmtip]{xy}
\DeclareFontFamily{U}{wncy}{}
\DeclareFontShape{U}{wncy}{m}{n}{<->wncyr10}{}
\DeclareSymbolFont{mcy}{U}{wncy}{m}{n}
\DeclareMathSymbol{\Sh}{\mathord}{mcy}{"58}
\usepackage[T2A,T1]{fontenc}
\usepackage[OT2,T1]{fontenc}

\usepackage{tikz}
\usetikzlibrary{shapes.geometric}
\usetikzlibrary{decorations.markings}
\tikzset{every loop/.style={min distance=10mm,looseness=10}}
\tikzstyle{vertex}=[auto=left,circle,minimum size=1pt,inner sep=0pt]

\newtheorem{theorem}{Theorem}[section]
\newtheorem{lemma}[theorem]{Lemma}

\newtheorem*{theorem*}{Theorem}
\newtheorem*{ass*}{Assumption}

\newtheorem{corollary}[theorem]{Corollary}
\newtheorem{remark}[theorem]{Remark}

\newtheorem{proposition}[theorem]{Proposition}

\newcommand{\Z}{\mathbb{Z}}

\newcommand{\Q}{\mathbb{Q}}
\newcommand{\F}{\mathbb{F}}

\newcommand{\rank}{\mathrm{rank}}

\DeclareFontEncoding{OT2}{}{}
\DeclareSymbolFont{cyrletters}{OT2}{wncyr}{m}{n}
\DeclareMathSymbol{\Sha}{\mathalpha}{cyrletters}{"58}

\numberwithin{equation}{section}

\begin{document}

\title{On Sufficient Conditions for Rank Zero in Families of Elliptic Curves}

\author[Shamik Das]{Shamik Das}
\address[Das]{Department of Mathematics and Statistics, IIT Kanpur, India}
\email{shamikd@iitk.ac.in}

\author[Debajyoti De]{Debajyoti De}
\address[De]{Department of Mathematics, IIT Madras, India}
\email{debajyotide20@gmail.com}

\author[Arup Mondal]{Arup Mondal}
\address[Mondal]{Department of Mathematics, IISER Pune, India}
\email{arup.mondal@students.iiserpune.ac.in}


\begin{abstract}
In this article, we study the Mordell–Weil ranks of elliptic curves in the family
$E_D:y^2=x^3+Dx$, where $D\in\mathbb{Z}\setminus\{0\}$. Under suitable congruence conditions on $D$ and its prime divisors, we use $2$-descent via a rational $2$-isogeny to study the rank of $E_D$.  Under our congruence hypotheses, the ranks of these matrices give sufficient conditions for $E_D$ to have rank zero. We give explicit sufficient conditions for the curves $E_D$ to have Mordell–Weil rank zero over $\mathbb{Q}$ in the following cases: (i) $D$ is an odd negative integer, (ii) $D$ is an odd positive integer, and (iii) $D$ is a positive even integer.  Also we give a quantitative result describing how many square-free parameters satisfy the required matrix-rank condition. In particular, our results show that these rank-zero criteria hold for a significant proportion of the admissible parameters.

\end{abstract}

\maketitle

\section{Introduction}

The study of rational points on elliptic curves has a long history in number theory. In the arithmetic of elliptic curves, the Mordell-Weil group plays a central role. By the Mordell-Weil theorem, the group of rational points $E(\Q)$ is finitely generated, isomorphic to $\Z^{r} \oplus E(\Q)_{\mathrm{tors}} $, where $r = \mathrm{rank}\,E(\Q)$ is the rank of the curve $E$ over $\Q$ and $E(\Q)_{\mathrm{tors}}$ is the torsion subgroup of $E(\Q)$. Since computing the rank $r$ of a general elliptic curve $E$ over $\mathbb{Q}$ is difficult, Selmer groups provide useful tools for bounding the rank.
In recent years, the study of families of elliptic curves of the form
$$
E_D: y^2=x^3+Dx,\qquad D\in\mathbb{Z}\setminus\{0\},
$$
has attracted significant attention. A notable example is obtained by taking $D=-n^2$, where $n$ is a positive square-free integer. In this case, $E_D$ is related to the classical congruent number elliptic curve.
A primary direction of research focuses on determining the Mordell–Weil rank of such curves for various forms of $D$. In 2007, Spearman \cite{Spearman} gave a sufficient condition
on an odd prime $p$ for the elliptic curve
$y^2=x^3-px$ to have rank two. In another paper \cite{Spearman2}, he gave a sufficient condition for the curve $y^2=x^3-2px$ to have rank three. Extending these investigations, Fujita and Terai \cite{Fujita-terai} analyzed curves of the form $y^2=x^3-p^k x$ for $k\in \{1,2,3\}$, establishing necessary and sufficient conditions for these curves to have rank one or two. Similarly, Daghigh and Didari \cite{Hasaan-Didari} determined the ranks of $y^2=x^3-3px$  and $y^2=x^3-pqx$, where $p$ and $q$ are distinct odd primes. Around the same time, Kim \cite{kim2015ranks} investigated the curves $y^2=x^3\pm4px$ for prime $p$. More recently, Mina and Bacani \cite{Mina-Bacani} examined $y^2=x^3-3pqx$ for distinct odd primes $p$ and $q$, establishing sufficient conditions for these curves to have rank zero and rank one, respectively. Recently, in \cite{Arkabrata-Bidisha-Richa} the authors used $2$-isogeny methods to obtain explicit possibilities for the Mordell–Weil ranks of $y^2=x^3-2pqx$ for distinct odd primes $p,q$, conditional on the Parity Conjecture.

For the curve $E_D$ several arithmetic criteria for determining the Mordell--Weil rank have been developed using Selmer groups, descent, and class groups. For the curve $E_D:y^2=x^3+Dx$, viewed as a quartic twist of $E_{1}:y^2=x^3+x$, \cite[Theorem 3.1 and 3.2]{DDM} shows that, for each odd prime $p\nmid D$, its $p$-Selmer rank is determined by its root number and a suitable character component of the relative $p$-class group of a cyclic quartic extension of a number field containing $\mathbb Q(i)$, whose isomorphism class is independent of $D$. Under the assumption $\omega(E_D)=+1$ considered here, their results provide computational criteria for establishing Mordell--Weil rank $0$ and finiteness of the Tate--Shafarevich group \cite[Corollary 3.3]{DDM}.
The celebrated conjecture of Birch and Swinnerton-Dyer (BSD) relates the Mordell--Weil rank $\operatorname{rank}E(\Q)$ of an elliptic curve $E/\Q$ to the order of vanishing of its $L$-function $L(E,s)$ at $s=1$. More precisely,
$\operatorname{ord}_{s=1}L(E,s)=\operatorname{rank}E(\Q).$ The Hasse--Weil conjecture predicts that $L(E,s)$ admits an analytic continuation to the whole complex plane and the completed $L$-function
$\Lambda(E,s)$ satisfies the functional equation $\Lambda(E,s)=\omega(E)\,\Lambda(E,2-s),$ where $\omega(E)\in \{\pm 1\}$ is called the global root number of $E$.
The global root number $\omega(E)$ is equivalently defined as the product of the local root numbers $\omega_v(E)\in\{\pm1\}$ over all places $v$ of $\Q$; namely,
$\omega(E)=\omega_\infty(E)\prod_{p}\omega_p(E),$
where $p$ runs over all finite primes. See \cite[Section 1]{ro} for a detailed discussion of global and local root numbers. The root number of $E_D: y^2=x^3+Dx$ can be calculated from the closed formula given in Birch-Stephens \cite[Equations. (9), and (10), Page 296]{BS}, $$\omega(E_D) = \omega_{\infty} \omega_2 \prod_{\substack{p \ge 3 \\ p^2 \mid \mid D}} \omega_p, $$
where $\omega_{\infty} = \operatorname{sgn}(D)$,\, $\omega_2 = \begin{cases}
    -1 & \text{ if } D \equiv 3, 5, 13, 15 \pmod{16},\\
    +1 & \text{ otherwise. }
\end{cases}$ \, and \, $\omega_p = \begin{cases}
    -1 & \text{ if } p \equiv 3 \pmod{4},\\
    +1 & \text{ otherwise}.
\end{cases}$

Now if we restrict ourselves to square-free integers $D$ then there is no prime $p$ for which
$p^2\parallel D$. Consequently, the product over odd primes equals $1$, since its indexing set is empty. Therefore,
$$\omega(E_D)=\omega_\infty\omega_2.$$
According to the sign and parity of $D$, this formula may be summarized as follows.

\begin{table}[h!]
\centering
\label{tab:root-numbers}
\begin{tabular}{l l c c l}
\toprule
\textbf{Sign of $D$} & \textbf{Parity of $D$} & $\omega_{\infty}$ & $\omega_2$ & \textbf{Root Number } $\omega(E_D)$ \\
\midrule
$D > 0$ & Odd & $+1$ & $\begin{cases} -1 & \text{if } D \equiv 3, 5, 13, 15 \pmod{16} \\ +1 & \text{if } D \equiv 1, 7, 9, 11 \pmod{16} \end{cases}$ & $\begin{cases} +1 & \text{if } D \equiv 1, 7, 9, 11 \pmod{16} \\ -1 & \text{if } D \equiv 3, 5, 13, 15 \pmod{16} \end{cases}$ \\
\addlinespace
$D > 0$ & Even & $+1$ & $+1$ & $+1$ \\
\addlinespace
$D < 0$ & Odd & $-1$ & $\begin{cases} -1 & \text{if } D \equiv 3, 5, 13, 15 \pmod{16} \\ +1 & \text{if } D \equiv 1, 7, 9, 11 \pmod{16} \end{cases}$ & $\begin{cases} +1 & \text{if } D \equiv 3, 5, 13, 15  \pmod{16} \\ -1 & \text{if } D \equiv 1, 7, 9, 11\pmod{16} \end{cases}$ \\
\addlinespace
$D < 0$ & Even & $-1$ & $+1$ & $-1$ \\
\bottomrule
\end{tabular}
\par\vspace{5pt}
\caption{Root numbers $\omega(E_D)$ according to the sign and parity of the square-free integer $D$.}
\end{table}
Note that, for an integer $D$ if the root number $\omega(E_D)=-1$, then according to the BSD conjecture, the Mordell--Weil rank of the elliptic curve $E_D$ over $\mathbb{Q}$ is expected to be odd and hence positive. Thus, it is natural to investigate the cases in which $\omega(E_D)=1$, where the root number does not predict whether the Mordell--Weil rank is zero or positive. In particular, it is interesting to determine for which families of integers $D$ the curve $E_D$ has Mordell--Weil rank zero. In this article, we investigate this question and obtain elliptic curves $E_D$ of Mordell--Weil rank zero for integers $D$ having arbitrarily many prime factors. More precisely, we give explicit sufficient conditions for the curves $E_D$ to have Mordell--Weil rank zero over $\mathbb{Q}$ in the following cases: (i) $D$ is an odd negative integer such that $-D\equiv \pm 3\pmod{8}$ (Theorem \ref{TheoremA negative D}); (ii) $D$ is an odd positive integer such that $D\equiv 7,11\pmod{16}$ (Theorem \ref{TheoremB positive D}); and (iii) $D$ is a positive even integer (Theorems \ref{thm-2D-C}, \ref{thm-2D-D} \& \ref{thm 2D-E}) satisfying certain restrictions on the prime factorization of $D$. We first introduce some notations that will be used in the statements of our main results and throughout the subsequent sections, including their proofs.

\textbf{Notations:} We will use the following terminologies. Let $r$ be a positive integer. 
\begin{itemize}
    \item  $\mathbf{0}_r=(0,\ldots,0)^{\mathsf T}$ and
         $\mathbf{1}_r=(1,\ldots,1)^{\mathsf T}$ denote the vectors consisting of $r$ zeros and $r$ ones, respectively.
    \item $M_r(\F_2)$ denotes the set of all $r\times r$ matrices over the field $\F_2$ where $\F_2$ is the finite field with two elements.
    \item $\mathcal I_r$ denotes the identity matrix of order $r$ over $\F_2$.
    \item For a vector $\boldsymbol{x}=(x_1,\ldots,x_r)^{\mathsf T}$, we denote its $i$-th entry by $(\boldsymbol{x})_i=x_i$, for $1\le i\le r$.
    \item For any matrix $A$ and any vector $\boldsymbol{v}$ of compatible dimensions, $A\cdot\boldsymbol{v}$ denotes their usual matrix-vector product.
    \item $\legendre{a}{p}$ is the Legendre symbol, where $p$ is an odd prime,  $a\in\Z$ and $p\nmid a$ .
    \item $\left[\frac{a}{p}\right]$ denotes the additive Legendre symbol, i.e., the image of $\legendre{a}{p}$ under the isomorphism as in \eqref{eq:phi-def}.
    
\end{itemize}

We now introduce the matrix that represents the quadratic residue information that is needed in our $2$-descent. For distinct odd primes $\ell_1, \ell_2, \ldots, \ell_r$, set $L := \ell_1 \ell_2 \cdots \ell_r$.  
    Define the matrix $\mathcal{A}_{L} = [a_{ij}] \in M_r(\F_2)$ by
    \begin{equation}\label{matrix A_D}
        a_{ij} =
        \begin{cases}
            1, & \text{if } \legendre{\ell_j}{\ell_i} = -1 \text{ and } i \neq j,\\[4pt]
            0, & \text{if } \legendre{\ell_j}{\ell_i} = 1 \text{ and } i \neq j,
        \end{cases}
        \qquad\text{and}\qquad
        a_{ii} = \sum_{j \ne i} a_{ij} .
    \end{equation}
    In the case $r = 1$, we set $\mathcal{A}_{L}$ to be the $1 \times 1$ zero matrix, $\mathcal{A}_{L} = [0]$.
    
We now state our main results concerning the Mordell--Weil rank of $E_D$. Under suitable congruence conditions on the prime divisors of $D$, we obtain sufficient conditions for $E_D(\Q)$ to have rank $0$. We also show that these families of elliptic curves have rank $0$ over the quadratic field $\Q(i)$ (See Corollary \ref{allCorollary}).  The curves under consideration are $E_{-D}$, $E_D$, or $E_{2D}$ according to sign or parity of $D$. Our results are as follows:

\begin{lthm}\label{TheoremA negative D}
Let $D = p_1 p_2 \cdots p_t q$ be a square-free integer, where $p_1,\ldots,p_t,q$ are primes and $t$ is a positive integer. Consider the elliptic curve $E_{-D}: y^2 = x^3 - Dx$ over $\Q$. Let $\mathcal{A}_{D} \in M_{t+1}(\F_2)$ be the matrix defined as in \eqref{matrix A_D}. Suppose one of the following conditions holds:
\begin{enumerate}[label=(\roman*)]
    \item $p_i \equiv 1 \pmod 8$ for all $i$, and $q \equiv 3 \pmod 8$;
    \item $p_i \equiv 5 \pmod 8$ for all $i$, $q \equiv 7 \pmod 8$, and $t$ is odd;
    \item $p_i \equiv 5 \pmod 8$ for all $i$, $q \equiv 3 \pmod 8$, and $t$ is even; 
    \item $p_i \equiv 1 \pmod 8$ for all $i$, $q \equiv 5 \pmod 8$, and $D \equiv 13 \pmod{16}$.
\end{enumerate}
If $\operatorname{rank}_{\mathbb{F}_2}(\mathcal{A}_{D}) = t$, then $\operatorname{rank} E_{-D}(\mathbb{Q}) = 0$.
\end{lthm}

\begin{lthm}\label{TheoremB positive D}
Let $D = p_1p_2\cdots p_tq$ be a square-free integer, where $p_1,\ldots,p_t,q$ are primes and $t$ is a positive integer. Suppose that $p_i \equiv 1 \pmod 4$ for all $i$, and $q \equiv 3 \pmod 4$, and $D \equiv 7 \text{ or } 11 \pmod{16}$. Consider the elliptic curve $E_D: y^2 = x^3 + Dx$ over $\Q$. Let $\mathcal{A}_{D} \in M_{t+1}(\F_2)$ be the matrix defined as in \eqref{matrix A_D}.  If $\operatorname{rank}_{\mathbb{F}_2}(\mathcal{A}_{D}) = t$, then $\operatorname{rank} E_D(\mathbb{Q}) = 0$.

\end{lthm}


\begin{lthm}
\label{thm-2D-C}
    Let $D = p_1 p_2 \cdots p_t q$ be a square-free integer, where $p_1,\ldots,p_t,q$ are primes and $t$ is a positive integer. Suppose that  $p_i \equiv 1 \pmod 8$ for $i \in \{1,2, \ldots, t\}$, and  $q \equiv \pm 3 \pmod 8$. Consider the elliptic curve $E_{2D}: y^2 = x^3 + 2Dx$ over $\Q$. Let $\mathcal{A}_D \in M_{t+1}(\F_2)$ be the matrix defined as in \eqref{matrix A_D}.  If $\operatorname{rank}_{\mathbb{F}_2}(\mathcal{A}_{D}) = t$, then $\operatorname{rank} E_{2D}(\mathbb{Q}) = 0$.
\end{lthm}

\begin{lthm}
\label{thm-2D-D}
    Let $D = p_1 p_2 \cdots p_t$ be a square-free integer, where $p_i$ are primes satisfying $p_i \equiv 5 \pmod{8}$ for all $i \in \{ 1, \ldots, t \}$ and $t$ is a positive integer. Consider the elliptic curve over $\Q$ given by $E_{2D}: y^2=x^3+2Dx.$ Let $\mathcal{A}_{D} \in M_{t}(\F_2)$ be the matrix defined as in \eqref{matrix A_D}.  If $\rank_{\F_2}(\mathcal{A}_{D}+\mathcal I_t)= t $, then $\operatorname{rank} E_{2D}(\Q) =0.$
\end{lthm}


\begin{lthm}\label{thm 2D-E}
Let 
$D=p_1p_2\cdots p_t$, where $p_1,\ldots,p_t$ are distinct primes satisfying $p_i\equiv3\pmod8$ for every $i$, and $t$ is a positive odd integer.
Consider the elliptic curve $E_{2D}:y^2=x^3+2Dx.$
Let $\mathcal A_D\in M_t(\F_2)$ be the matrix defined as in \eqref{matrix A_D}.
If $\operatorname{rank}_{\F_2}(\mathcal A_D)=t-1,$
then $\operatorname{rank}E_{2D}(\Q)=0$.
\end{lthm}

\begin{remark}
The conditions in Theorems \ref{TheoremA negative D}--\ref{thm 2D-E} are sufficient, but not necessary, for the corresponding elliptic curves to have Mordell--Weil rank $0$ over $\Q$. To illustrate this, we consider the following examples.

\begin{enumerate}[label=(\roman*)]

\item[Theorem \ref{TheoremA negative D}:]  Consider $D=17\cdot 19.$
In this case,  $\operatorname{rank}_{\F_2}(\mathcal{A}_{D})=0<1, \text{ but } \operatorname{rank}E_{-D}(\Q)=0.$

\item[Theorem \ref{TheoremB positive D}:] Consider $D=41\cdot 31.$
In this case,  $\operatorname{rank}_{\F_2}(\mathcal{A}_{D})=0<1, \text{ but } \operatorname{rank}E_{D}(\Q)=0.$ 

\item[Theorem \ref{thm-2D-C}:]
Consider $D=17 \cdot  41\cdot 43.$
In this case,  $\operatorname{rank}_{\F_2}(\mathcal{A}_{D})=1<2, \text{ but } \operatorname{rank}E_{2D}(\Q)=0.$

\item[Theorem \ref{thm-2D-D}:] Consider $D=5\cdot13 \cdot 149.$ In this case,  $\operatorname{rank}_{\F_2}(\mathcal{A}_{D}+\mathcal I_t)=2<3, \text{ but } \operatorname{rank}E_{2D}(\Q)=0.$

\item[Theorem \ref{thm 2D-E}:] 
Consider $D=3\cdot 11\cdot 67.$ In this case,  $\operatorname{rank}_{\F_2}(\mathcal{A}_{D})=1<2, \text{ but } \operatorname{rank}E_{2D}(\Q)=0.$

\end{enumerate} 

\end{remark}

\textbf{Motivation and Organisation:}  The main results are proved using explicit $2$-descent via a rational $2$-isogeny. We determine the images of the corresponding descent maps by studying the associated quartic Diophantine equations. For square-free $D$, the local conditions arising in this descent are controlled by quadratic residue relations among the prime divisors of $D$. These relations can be expressed in a matrix over $\F_2$ using the Legendre symbols. Thus the rank of this matrix provides a linear algebra framework for studying the Mordell--Weil rank of $E_D$.
The main results (Theorems \ref{TheoremA negative D}--\ref{thm 2D-E}) are proved using explicit $2$-descent and associated quartic equations. The paper is organised as follows. In Section~\ref{Prelimary section}, we recall the necessary preliminaries on $2$-descent, and establish the matrix-theoretic framework used in the subsequent section. The proofs of the main theorems are given in Section~\ref{Proof of main thms}. In Section~\ref{Density section}, we prove the corresponding density results for the square-free parameters satisfying the required matrix-rank conditions. Finally, Section~\ref{Numerical Section} presents numerical data illustrating the theoretical results.


\section{Preliminaries}\label{Prelimary section}

In this section, we recall some fundamental results and discuss the tools of $2$-descent via $2$-isogeny, which will be used to compute the ranks of the elliptic curves over $\Q$. We also establish some preliminary results in linear algebra that play a crucial role in the proof of our main theorem.

\subsection{\texorpdfstring{$2$-Descent via $2$-Isogeny}{lg}}

Suppose that $E : Y^2 = X^3 + aX^2 + bX$ is an elliptic curve over $\Q$ where $b\neq 0$ and $a^2-4b\neq 0$. Let $E' : Y^2 = X^3 - 2aX^2 + (a^2 - 4b)X$ be the corresponding $2$-isogenous elliptic curve to $E$ and the corresponding isogeny map for $x\neq 0$ is given by 
\begin{equation}\label{2 isogeny}
    \Phi : E \longrightarrow E', \qquad \Phi(x,y)=\left(\frac{y^2}{x^2},\, \frac{y(x^2-b)}{x^2}\right)
\end{equation}
 which is of degree $2$.
Let $\Q^{*}$ denote the multiplicative group of nonzero rational numbers,
and let $\Q^{*2}$ be its subgroup consisting of rational squares.
Thus, $\Q^{*}/\Q^{*2}$ denotes the group of nonzero rational numbers
modulo squares. We write $E(\Q)$ and $E'(\Q)$ for the groups of
rational points on $E$ and its $2$-isogenous curve $E'$, respectively. We
define the $2$-descent homomorphism $\alpha:E(\mathbb{Q})\longrightarrow \mathbb{Q}^{*}/\mathbb{Q}^{*2}$ by
\begin{equation}
\label{alpha}
 \alpha(P)=
\begin{cases}
1 \pmod{\mathbb{Q}^{*2}}, & \text{if } P=\mathcal{O}_E, \text{ the point at infinity},\\[4pt]
b \pmod{\mathbb{Q}^{*2}}, & \text{if } P=(0,0),\\[4pt]
x \pmod{\mathbb{Q}^{*2}}, & \text{if } P=(x,y)\text{ with }x\neq 0.
\end{cases}   
\end{equation}
Similarly, we define the $2$-descent homomorphism on the $2$-isogenous
curve $E'$ by $
\alpha':E'(\Q)\longrightarrow
\Q^{*}/\Q^{*2}$
as follows:
\begin{equation}
\label{alpha'}
\alpha'(P)=
\begin{cases}
1 \pmod{\Q^{*2}}, & \text{if } P=\mathcal{O}_{E'}, \text{ the point at infinity},\\[4pt]
b' \pmod{\Q^{*2}}, & \text{if } P=(0,0),\\[4pt]
x \pmod{\Q^{*2}}, & \text{if } P=(x,y)\text{ with }x\neq 0,
\end{cases}
\end{equation}
where $b'= a^2-4b$ is the coefficient of $X$ in the equation of $E'$.

The descent maps \eqref{alpha} and \eqref{alpha'} reduce the problem
of determining the rank of an elliptic curve over $\Q$ to the study
of their images in $\Q^*/\Q^{*2}$. The following standard results
provide a method for computing the rank (see
\cite{silverman2009arithmetic} for details).

\begin{proposition}\label{rank proposition}
Let $E/\Q$ be an elliptic curve admitting a $2$-isogeny $\Phi:E\longrightarrow E'$ defined over $\Q$ as in \eqref{2 isogeny}. Then the Mordell--Weil rank of $E(\Q)$ satisfies
$$2^{\operatorname{rank} E(\Q)+2}
=\left|\operatorname{Im}(\alpha)\right|
\cdot\left|\operatorname{Im}(\alpha^\prime)\right|,$$
where $\alpha$ and
$\alpha^\prime$
are the descent maps defined in \eqref{alpha} and
\eqref{alpha'} respectively.
\end{proposition}

By Proposition \ref{rank proposition}, the rank of an elliptic curve
$E$ of the given form can be determined by computing the cardinalities
of the images of the descent maps \eqref{alpha} and \eqref{alpha'}. The following result provides
a criterion for determining whether an element belongs to these
images. For a proof, see \cite{silverman2009arithmetic}.

\begin{lemma}
\label{image-torsors}
Let $E:Y^2=X^3+aX^2+bX$ and $E':Y^2=X^3-2aX^2+(a^2-4b)X$ be the corresponding $2$-isogenous elliptic curves, where
$a,b\in\mathbb Z$ and $b(a^2-4b)\neq0$.
Then the following hold:
\begin{itemize}
    \item Every element of $\operatorname{Im}(\alpha)$
    is represented by a square-free divisor of $b$.
    For each such divisor $d$, we have
    \[
    d\in\operatorname{Im}(\alpha)
    \quad\Longleftrightarrow\quad
    N^2=dM^4+aM^2e^2+\frac{b}{d}e^4
    \]
    for some $(N,M,e)\in\mathbb Z^3$ satisfying
    $\gcd(M,e)=1$.

    \item Every element of $\operatorname{Im}(\alpha')$
    is represented by a square-free divisor of $a^2-4b$.
    For each such divisor $d$, we have
    \[
    d\in\operatorname{Im}(\alpha')
    \quad\Longleftrightarrow\quad
    N^2=dM^4-2aM^2e^2+
    \frac{a^2-4b}{d}e^4
    \]
    for some $(N,M,e)\in\mathbb Z^3$ satisfying
   $\gcd(M,e)=1$.
\end{itemize}
\end{lemma}

Consider the elliptic curve $E_D:y^2=x^3+Dx$, where
$D\in\mathbb Z\setminus\{0\}$. To determine the rank of $E_D$, we study the images of the descent maps $\alpha$ and $\alpha'$, defined in \eqref{alpha} and \eqref{alpha'}  respectively. By Lemma \ref{image-torsors}, these images can be determined by studying the solvability of the corresponding quartic Diophantine equations in integers $N,M,e$ with $\gcd(M,e)=1$. For $E_{-D}$ with $D>0$ square-free, the following lemma
shows that it suffices to determine the positive square classes in $\operatorname{Im}(\alpha)$.

\begin{lemma}
\label{E{-D}-d-equivalence}

 Let $D>0$ be a square-free integer, and let $E_{-D}:y^2=x^3-Dx$. Consider the $2$-descent map $\alpha:E_{-D}(\Q)\longrightarrow \Q^{*}/\Q^{*2}$ as in \eqref{alpha}. Then, for every positive divisor $d\mid D$,
$$
\pm d\in\operatorname{Im}(\alpha) \quad\Longleftrightarrow\quad \mp\frac{D}{d}\in\operatorname{Im}(\alpha).
$$
In particular, the negative elements in $\operatorname{Im}(\alpha)$ are uniquely determined by the positive elements and conversely.
\end{lemma}

\begin{proof}
   Since $H=\operatorname{Im}(\alpha)$ is a subgroup of $\Q^*/\Q^{*2}$ containing $\alpha((0,0))=-D$, for every positive divisor $d\mid D$, we have\
$\pm d\in H \iff (-D)(\pm d)\in H \iff \mp\frac{D}{d}\in H.$
This proves the result.
\end{proof}

The previous lemma shows that the negative square classes in the image of the $2$-descent map for $E_{-D}$ are determined by the positive ones in image of $\alpha$, and vice versa.  We now introduce two group isomorphisms that identify signs and square-free divisors with elements and vectors over $\F_2$, respectively, thereby translating multiplicative relations into additive ones. Consider the group isomorphism
\begin{equation}\label{eq:phi-def}
\phi:\{\pm 1\} \longrightarrow \{0,1\}, \quad \text{ defined by } \phi\bigl((-1)^{\epsilon}\bigr)=\epsilon,
\quad \text{for } \epsilon \in \{0,1\}.
\end{equation} 
In this context, $\{\pm 1\}$ is regarded as a multiplicative group, while $\{0,1\}$ is viewed as an additive group modulo $2$.  

Let $L = \ell_{1} \cdots \ell_{r}$ be a positive square-free integer, where $\ell_{1},\dots,\ell_{r}$ are distinct odd primes. Consider the set $G(L)$ of all positive divisors of $L$ i.e. $G(L) =\{a: ~a> 0, ~a\mid L\}.$ A straightforward calculation yields that $G(L)$ has a group structure under the following binary operation
\begin{equation*}\label{thegroupG}
u*u' := \frac{uu'}{(\gcd(u,u'))^2}, \quad  u, u^\prime \in G(L).
\end{equation*}
Every $u \in G(L)$ can be written uniquely in the form
$u = \prod_{i=1}^{r} \ell_{i}^{\varepsilon_{i}(u)},$
where each exponent $\varepsilon_{i}(u) \in \{0,1\}$. Define a map 
\begin{align}
\label{varphi-def} 
\varphi : G(L) \longrightarrow \F_2^{r}  ~\text{ by }~  \varphi(u)= \big(\varepsilon_{1}(u), \dots, \varepsilon_{r}(u)\big)^{\mathsf T},  
\end{align}
 where $\varepsilon_{i}(u)=1$ if $\ell_{i}\mid u$ and  $\varepsilon_{i}(u)=0$  otherwise. One checks easily that $\varphi$ is a group isomorphism.  We next discuss some elementary linear algebra results over $\F_2$ that will be used in the proofs of the main theorems. We begin with a characterization of the image of a symmetric matrix whose row sums are zero. 

\begin{lemma}\label{Image of A_D}
Let $\mathcal{A}\in M_{r}(\F_2)$ be a symmetric matrix. Suppose that $\operatorname{rank}_{\F_2}(\mathcal{A})=r-1$ and $\mathcal{A}\cdot \mathbf{1}_{r}=\mathbf{0}_{r}.$ Then $\operatorname{Im}(\mathcal{A})=\left\{(x_1,\ldots,x_{r})^{\mathsf T}
\in\mathbb{F}_2^{r}:\displaystyle\sum_{i=1}^{r}x_i=0\right\}.$ Thus, $\operatorname{Im}(\mathcal{A})$ consists precisely of those vectors having an even number of nonzero entries.
\end{lemma}

 \begin{proof}
 Since $\operatorname{rank}_{\F_2}(\mathcal{A})=r-1$, we have $\dim\ker(\mathcal{A})=1.$
 As $\mathbf{1}_{r}\in\ker(\mathcal{A})$, we have $\ker(\mathcal{A})=\langle\mathbf{1}_{r}\rangle.$ Since $\mathcal{A}$ is symmetric, $ \operatorname{Im}(\mathcal{A})=\ker(\mathcal{A})^\perp,$ which gives the description of $\operatorname{Im}(\mathcal{A})$ as required.
 \end{proof}

\begin{lemma}\label{rowsumzero for p_i}
Let $r\ge1$, and let $L=\ell_1\cdots\ell_r$ be a product of distinct odd primes. Let $\mathcal A_L=(a_{ij})$ be the matrix defined in \eqref{matrix A_D}. For $\ell\in G(L)$, write $\ell=\ell_1^{e_1}\cdots\ell_r^{e_r}, \;\; e_i\in\{0,1\},$ and associate to $\ell$ the vector $\boldsymbol v_\ell =\varphi(\ell)
  =(e_1,\ldots,e_r)^{\mathsf T}\in\F_2^r$ where $\varphi$ is the group isomorphism defined in
\eqref{varphi-def}. Let $\lambda,\mu\in\Z\setminus\{0\}$ satisfy
$\gcd(\lambda\mu,L)=1$. Suppose that the equation
$$C_\ell:\quad N^2=\lambda\ell M^4+\mu\frac{L}{\ell}e^4$$
has an integer solution $(N,M,e)$ with $\gcd(M,e)=1$. Then, for every $1\le i\le r$,
$$(\mathcal A_L\cdot\boldsymbol v_\ell)_i
=\begin{cases}
\left[\dfrac{\mu}{\ell_i}\right],
    & \text{if }\ell_i\mid\ell,\\[6pt]
\left[\dfrac{\lambda}{\ell_i}\right],
    & \text{if }\ell_i\nmid\ell.
\end{cases}
$$
\end{lemma}

\begin{proof}
    All matrix operations and sums of additive Legendre symbols are taken in $\F_2$. Let $(N,M,e)$ be an integer solution with $\gcd(M,e)=1$, of
$N^2=\lambda\ell M^4+\mu\frac{L}{\ell}e^4.$ Fix $i\in\{1,\ldots,r\}$.  By the definition of $\mathcal A_L$  in \eqref{matrix A_D} we have,
\begin{equation}
\label{matrix-identity}
 (\mathcal A_L \cdot \boldsymbol v_\ell)_i =a_{ii}e_i+\sum_{j\ne i}a_{ij}e_j
=\sum_{j\ne i}a_{ij}(e_i+e_j)
=\sum_{j\ne i} \left[\frac{\ell_j}{\ell_{i}}\right](e_i+e_j).   
\end{equation}

First suppose that $\ell_{i}\mid\ell$, so that $e_i=1$.
The preceding identity \eqref{matrix-identity} gives
$$
(\mathcal A_L\cdot \boldsymbol v_\ell)_i
=\sum_{j\ne i}
  \left[\frac{\ell_j}{\ell_{i}}\right](1+e_j)
  =\sum_{\substack{j\ne i\\ e_j=0}}
  \left[\frac{\ell_j}{\ell_{i}}\right]
=\left[\frac{L/\ell}{\ell_{i}}\right].
$$
The last equality follows since $L/\ell$ is the product of precisely those primes $\ell_j$ for which $e_j=0$. We claim that $\ell_i\nmid e$. Indeed, if $\ell_i\mid e$, then $\ell_i\nmid M$ because $\gcd(M,e)=1$. Since $\ell$ is square-free and $\ell_i\nmid\lambda$, the term $\lambda\ell M^4$ is divisible by $\ell_i$ but not by $\ell_i^2$, whereas $\mu(L/\ell)e^4$ is divisible by $\ell_i^4$. Thus the right-hand side of $C_\ell$ is divisible by $\ell_i$ but not by $\ell_i^2$, which is impossible. Hence $\ell_i\nmid e$. Reducing the equation $C_{\ell}$ modulo $\ell_{i}$, we obtain
$$N^2\equiv \mu \frac{L}{\ell}e^4\pmod{\ell_{i}} \implies \legendre{\mu(L/\ell)}{\ell_{i}}=1.$$
 Consequently, we have $\left[\frac{\mu}{\ell_{i}}\right] + \left[\frac{L/\ell}{\ell_{i}}\right]=0$. Therefore, $$(\mathcal A_L \cdot\boldsymbol v_\ell)_i
=\left[\frac{L/\ell}{\ell_{i}}\right]=\left[\frac{\mu}{\ell_{i}}\right].$$
Now suppose that $\ell_{i}\nmid\ell$, so that $e_i=0$. The identity \eqref{matrix-identity}  gives
$$(\mathcal A_L\cdot\boldsymbol v_\ell)_i
=\sum_{j\ne i}\left[\frac{\ell_j}{\ell_{i}}\right]e_j
=\left[\frac{\ell}{\ell_{i}}\right].$$
 We claim that $\ell_i\nmid M$.
If $\ell_i\mid M$, then $\ell_i\nmid e$
because $\gcd(M,e)=1$. Since $L/\ell$ is square-free and $\ell_i\nmid\mu$, the term $\mu(L/\ell)e^4$ is divisible by $\ell_i$
but not by $\ell_i^2$, while $\lambda\ell M^4$
is divisible by $\ell_i^4$. Again, the right-hand side of $C_\ell$ cannot be a square. Thus $\ell_i\nmid M$. 

Reducing the defining equation modulo $\ell_{i}$, we obtain
$$N^2\equiv \lambda\ell M^4\pmod{\ell_{i}} \implies \legendre{\lambda\ell}{\ell_{i}}=1.$$
Equivalently,
$\left[\frac{\lambda}{\ell_{i}}\right]+ \left[\frac{\ell}{\ell_{i}}\right]=0.$
Hence $(\mathcal A_L\cdot\boldsymbol v_\ell)_i= \left[\frac{\lambda}{\ell_{i}}\right]$.
We are done.
\end{proof}

\begin{lemma}\label{odd lemma invertible}
 Let $r$ be a positive odd integer, and let
$\mathcal{A}=(a_{ij})\in M_r(\F_2)$ satisfy
$$a_{ij}+a_{ji}=1 \quad (i\ne j),
\quad a_{ii}=\sum_{j\ne i}a_{ij}.$$
If $\operatorname{rank}_{\F_2}(\mathcal{A})=r-1$, then
$\mathcal{A}+\mathcal I_r$ is invertible.
\end{lemma}

\begin{proof}
    Let $J=\mathbf{1}_r \cdot \mathbf{1}_r^{\mathsf T}$ be the $r \times r$ all-ones matrix. Since $a_{ij}+a_{ji}=1$ for $i\ne j$, we have $\mathcal{A}+\mathcal{A}^{\mathsf T}=\mathcal I_r+J$. 

    We first show that $\mathcal{A} \cdot \mathbf{1}_r=\mathcal{A}^{\mathsf T}\cdot \mathbf{1}_r=\mathbf{0}_r$. By the definition of the diagonal entries of $\mathcal{A}$, the sum of the entries in the $i$th row is $a_{ii}+\sum_{j\neq i}a_{ij}=2a_{ii}=0.$ Thus, $\mathcal{A}\cdot \mathbf{1}_r=\mathbf{0}_r$. On the other hand, since $a_{ij}=1+a_{ji}$ for $i\ne j$
    and $r$ is odd, we obtain $\sum_{i\ne j}a_{ij}
    =(r-1)+\sum_{i\ne j}a_{ji} =a_{jj}.$ Consequently, the sum of the entries in the $j$th column is $2a_{jj}=0$, yielding
   $\mathcal A^{\mathsf T}\cdot \mathbf{1}_r=\mathbf{0}_r$.

   Since $\operatorname{rank}\mathcal{A}=\operatorname{rank}\mathcal{A}^{\mathsf T}=r-1$, both kernels are one-dimensional, and therefore $\ker \mathcal{A}=\ker \mathcal{A}^{\mathsf T}=\operatorname{span}\{\mathbf{1}_r\}.$
   Now let $\boldsymbol{x}\in\mathbb{F}_2^r$ satisfy
  $(\mathcal A+\mathcal I_r)\cdot\boldsymbol{x}=\mathbf{0}_r$.
 Then $\mathcal A \cdot \boldsymbol{x}=\boldsymbol{x}$, and hence
  $\mathbf{1}_r^{\mathsf T}\cdot\boldsymbol{x}
  =\mathbf{1}_r^{\mathsf T}\cdot \mathcal A\cdot\boldsymbol{x}
  =0,$ because $\mathbf{1}_r^{\mathsf T}\cdot \mathcal{A}
  =(\mathcal{A}^{\mathsf T}\cdot\mathbf{1}_r)^{\mathsf T}
   =\mathbf{0}_r^{\mathsf T}.$

    Hence $J\cdot \boldsymbol{x}=\mathbf{0}_r$, and so $$\mathcal{A} \cdot \boldsymbol{x}+\mathcal{A}^{\mathsf T}\cdot \boldsymbol{x}=(\mathcal{A}+\mathcal{A}^{\mathsf T}) \cdot\boldsymbol{x}=(\mathcal I_r+J)\cdot\boldsymbol{x}=\boldsymbol{x}.$$
    Since $\mathcal{A} \cdot \boldsymbol{x}=\boldsymbol{x}$, this forces $\mathcal{A}^{\mathsf T}\cdot \boldsymbol{x}= \mathbf{0}_r$, so $\boldsymbol{x}\in\ker \mathcal A^{\mathsf T}=\operatorname{span}\{\mathbf{1}_r\}$. Thus $\boldsymbol{x}=\mathbf{0}_r$ or $\boldsymbol{x}=\mathbf{1}_r$. As $r$ is odd, we have $\mathbf{1}_r^{\mathsf T}\cdot \mathbf{1}_r=1$
   in $\mathbb{F}_2$. Hence the second possibility contradicts $\mathbf{1}_r^{\mathsf T} \cdot \boldsymbol{x}=0$. Hence $\boldsymbol{x}=\mathbf{0}_r$. Therefore $\ker(\mathcal{A}+\mathcal I_r)=\{\mathbf{0}_r\}$, and by the rank--nullity theorem $\mathcal{A}+\mathcal I_r$ is invertible.
\end{proof}

\section{Proofs of Main Results}\label{Proof of main thms}
In this section, we prove the main theorems using the results established in the preliminary section. The proofs combine the descent criteria with the properties of the associated matrices to obtain the desired rank statements.

\begin{proof}[\textbf{\textit{Proof of Theorem \ref{TheoremA negative D}}}]
Let $E=E_{-D}:y^2=x^3-Dx$ and let $E^\prime=E_{4D}:y^2=x^3+4Dx$ be the curve $2$-isogenous to $E$. Let $\alpha$ and $\alpha^\prime$ be the corresponding descent maps as defined in \eqref{alpha} and \eqref{alpha'} respectively. We prove that
$$\operatorname{Im}(\alpha)=\{1,-D\},\quad \operatorname{Im}(\alpha^\prime)=\{1,D\}$$ and apply Proposition \ref{rank proposition} to get our desired result.

Every element of $\operatorname{Im}(\alpha)$ admits a representative of the form $\pm d$, where $d$ is a positive divisor of $D$. By Lemma \ref{E{-D}-d-equivalence}, it  suffices to determine the positive divisors $d$ of $D$ that belong to $\operatorname{Im}(\alpha)$.
 Let $d$ be a positive divisor of $D$ representing a class
in $\operatorname{Im}(\alpha)$.
Write $d=p_1^{e_1}\cdots p_t^{e_t}q^f,$ where $ e_1,\ldots,e_t,f\in\{0,1\},$
and set $\boldsymbol v_d=\varphi(d)
=(e_1,\ldots,e_t,f)^{\mathsf T}\in\F_2^{t+1},$
where $\varphi$ is the group isomorphism defined in \eqref{varphi-def}.
 Since $\operatorname{rank}_{\F_2}(\mathcal{A}_{D})=t$ and the sum of every row of $\mathcal{A}_{D}$ is zero, we have that $\ker(\mathcal{A}_{D})=\langle\mathbf{1}_{t+1} \rangle$.
By Lemma~\ref{image-torsors}, the associated quartic equation
\begin{equation}\label{Torsor for thm A}
C_d:N^2=dM^4-\frac{D}{d}e^4,
\end{equation}
has an integer solution $(N,M,e)$ with $\gcd(M,e)=1$.
In every case of the theorem, $p_i\equiv1\pmod4$
for $1\le i\le t$. Hence
$\left[\frac{-1}{p_i}\right]=0$
and $\left[\frac{1}{p_i}\right]=0.$
Applying Lemma~\ref{rowsumzero for p_i} with
$(\lambda,\mu)=(1,-1)$, we obtain
\begin{equation}\label{ith row=0}
    (\mathcal A_D\cdot\boldsymbol v_d)_i=0,
\qquad 1\le i\le t.
\end{equation}

First, we show that, in each of the four cases of Theorem \ref{TheoremA negative D}, the existence of integer solution $(N,M,e)$ with $\gcd(M,e)=1$ to $C_d$, as defined in \eqref{Torsor for thm A}, implies that $d=1$.

\textbf{(i)} Suppose that $p_i\equiv1\pmod8$ for all $i \in \{1, \ldots,t\}$ and $q\equiv3\pmod8$. 

Suppose that $q\mid d$. Then $f=1$, and hence $D/d\equiv1\pmod8$. From \eqref{Torsor for thm A}, we see that exactly one of $N$, $M$, and $e$ is even, while the other two are odd. If $e$ is even, then reducing \eqref{Torsor for thm A} modulo $8$ gives $N^2\equiv3\pmod8$, which is impossible. The other two cases can be ruled out similarly. Therefore, $q\nmid d$, and hence $f=0$. Reducing \eqref{Torsor for thm A} modulo $q$ gives $\legendre{d}{q}=1,$ and hence $(\mathcal{A}_{D}\cdot\boldsymbol{v}_{d})_{t+1}=0$. Thus $\mathcal{A}_{D}\cdot\boldsymbol{v}_{d}=\mathbf{0}_{t+1}$, so $\boldsymbol{v}_{d}\in\ker(\mathcal{A}_{D})$. Since $f=0$ and $\mathcal{A}_D$ has nullity $1$, we obtain $\boldsymbol{v}_d=\mathbf{0}_{t+1}$, which implies $d=1$.

\textbf{(ii)} Suppose that $p_i\equiv5\pmod8$ for all $i \in \{1, \ldots,t\}$, $q\equiv7\pmod8$, and $t$ is odd. It remains to show that $(\mathcal{A}_{D}\cdot\boldsymbol{v}_{d})_{t+1}=0.$

Suppose that $q\mid d$, then by Lemma \ref{rowsumzero for p_i} we have that $(\mathcal{A}_{D}\cdot\boldsymbol{v}_{d})_{t+1}=1$ as $(-1/q)=-1$.
Together with \eqref{ith row=0} this gives $\mathcal{A}_{D}\cdot\boldsymbol{v}_{d}=(0,\ldots,0,1)^{\mathsf T}$;  we get that $\mathbf{1}_{t+1}^{\mathsf T}\cdot
(\mathcal{A}_{D}\cdot\boldsymbol{v}_{d})=1$. Since $p_i \equiv 1 \pmod 4$ for every $i$, quadratic reciprocity implies that $\mathcal A_D$ is symmetric. Consequently,  $\mathbf{1}_{t+1}^{\mathsf T}\cdot \mathcal{A}_{D}=\mathbf{0}_{t+1}^{\mathsf T}$. This is a contradiction.  Hence $q\nmid d$, so $f=0$. Therefore, we obtain $(\mathcal{A}_{D}\cdot\boldsymbol{v}_{d})_{t+1}=0$. Thus $\boldsymbol{v}_{d}\in\ker(\mathcal{A}_{D})$, and $f=0$ gives $\boldsymbol{v}_{d}=\mathbf{0}_{t+1}$. Hence $d=1$.

\textbf{(iii)} Suppose that $p_i\equiv5\pmod8$ for all $1\le i\le t$,
$q\equiv3\pmod8$, and $t$ is even.
By an argument similar to that in part (ii), we obtain $d=1$.

\textbf{(iv)} Suppose that $p_i\equiv1\pmod8$ for all $i \in \{1, \ldots,t \}\;\;$, $q\equiv5\pmod8$, and $D\equiv13\pmod{16}$. By Lemma \ref{rowsumzero for p_i}, we obtain  $(\mathcal{A}_{D}\cdot\boldsymbol{v}_{d})_{i}=0$ for each $1 \le i \le t+1$. Since $\ker(\mathcal{A}_{D})=\langle \mathbf{1}_{t+1} \rangle,$ either $d=1$ or $d=D$. 

If $d=D$, then \eqref{Torsor for thm A} reduces to $N^2=DM^4-e^4.$ If $M,e$ are both odd, then $N^2\equiv D-1\equiv12\pmod{16}$; if $M$ is even and $e$ is odd, then $N^2\equiv-1\pmod{16}$; and if $M$ is odd and $e$ is even, then $N^2\equiv D\equiv13\pmod{16}$. Each case is impossible, since the only quadratic residues modulo $16$ are $0,1,4$ and $9$. Thus $d=1$.
Hence by Lemma \ref{E{-D}-d-equivalence} we obtain,
$\operatorname{Im}(\alpha)=\{1,-D\}.$

We now consider the isogenous curve $E^\prime$. By Lemma~\ref{image-torsors}, the associated quartic equation
\begin{equation}\label{dual torsor}
C_{d'}':N^2=d'M^4+\frac{4D}{d'}e^4.
\end{equation}

 has an integer solution.  Since $D>0$, it follows readily that $d'>0$. For a positive odd divisor $d^\prime \mid D$, write $d^\prime =p_1^{e_1}\cdots p_t^{e_t}q^f,$ and  $\boldsymbol{v}_{d^\prime}=(e_1,\ldots,e_t,f)^{\mathsf T}$ as in \eqref{varphi-def}. Lemma~\ref{rowsumzero for p_i}, applied with
$(\lambda,\mu)=(1,4)$, gives $\mathcal A_D\cdot \boldsymbol v_{d'}=\mathbf0_{t+1}.$ Hence $\boldsymbol v_{d'}$ is either $\mathbf0_{t+1}$
or $\mathbf1_{t+1}$, so $d'=1$ or $d'=D$. It remains to consider the even classes $2d'$ with $d'\mid D$. In this case \eqref{dual torsor} becomes
$N^2=2d'M^4+\frac{2D}{d'}e^4.$
By applying Lemma~\ref{rowsumzero for p_i} we obtain,
$$\mathcal{A}_{D}\cdot\boldsymbol{v}_{d^\prime}= \boldsymbol{b}, \quad
\boldsymbol b=
\left(
\left[\frac{2}{p_1}\right],
\ldots,
\left[\frac{2}{p_t}\right],
\left[\frac{2}{q}\right]
\right)^{\mathsf T}.
$$

By Lemma~\ref{Image of A_D}, every vector in $\operatorname{Im}(\mathcal{A}_{D})$ has coordinate sum zero in $\F_2$. In cases (i) and (iv), $\boldsymbol{b}=(0,\ldots,0,1)^{\mathsf T}$, whose coordinate sum is $1$. In case (ii), $\boldsymbol{b}=(1,\ldots,1,0)^{\mathsf T}$; since $t$ is odd, its coordinate sum is $t \equiv 1 \pmod2$. In case (iii), $\boldsymbol{b}=(1,\ldots,1)^{\mathsf T}$; since $t$ is even, its coordinate sum is $t+1 \equiv 1 \pmod2$. Thus $\boldsymbol{b}\notin\operatorname{Im}(\mathcal{A}_{D})$ in all four cases.
Thus no even class belongs to $\operatorname{Im}(\alpha')$.
Since every class in this image is positive and its odd classes are precisely $1$ and $D$, we conclude that,
$\operatorname{Im}(\alpha^\prime)=\{1,D\}.$ Finally, Proposition~\ref{rank proposition} gives 
$$2^{\operatorname{rank}E_{-D}(\Q)+2}
=|\operatorname{Im}(\alpha)| \,|\operatorname{Im}(\alpha')| =4.$$
Consequently, $\operatorname{rank}E_{-D}(\Q)=0$. This proves the theorem.
\end{proof}

\begin{proof}[\textbf{\textit{Proof of Theorem \ref{TheoremB positive D}}}]
Let $E=E_D:y^2=x^3+Dx$ and the elliptic curve $2$-isogenous to $E$ is $E^\prime=E_{-4D}:y^2=x^3-4Dx$. Let
$\alpha, \alpha^\prime$ be the descent maps defined in \eqref{alpha} and \eqref{alpha'} respectively. We prove that
$\operatorname{Im}(\alpha)=\{1,D\}, \;\;
\operatorname{Im}(\alpha^\prime)=\{1,-D\}.$ 

We first determine $\operatorname{Im}(\alpha)$. Since $\operatorname{rank}_{\F_2}(\mathcal{A}_{D})=t$ and every row of $\mathcal{A}_{D}$ has sum zero,
$\ker \mathcal{A}_{D}=\langle\mathbf{1}_{t+1} \rangle.$ By Lemma \ref{image-torsors}, every class in $\operatorname{Im}(\alpha)$ is represented by a  square-free divisor $d$ of $D$ for which $C_d:~N^2=dM^4+\frac{D}{d}e^4$ has an integer solution $(N,M,e)$ with $\gcd(M,e)=1$. If $d<0$, then $\frac{D}{d}<0$, as $D>0$, and the right-hand side of $C_{d}$ is strictly negative because $M$ and $e$ are not both zero. Thus $d>0$. For a positive divisor  $d=p_1^{e_1}\cdots p_t^{e_t}q^f$ of $D$, put $\boldsymbol{v}_d=(e_1,\ldots,e_t,f)^{\mathsf T}\in\mathbb F_2^{t+1}$ as in \eqref{varphi-def}. By Lemma~\ref{rowsumzero for p_i} with
$(\lambda,\mu)=(1,1)$ we obtain that $\mathcal A_D\cdot \boldsymbol{v}_d=\mathbf{0}_{t+1}$. Since
$\ker(\mathcal A_D)=\{\mathbf{0}_{t+1},\mathbf{1}_{t+1}\}$,
we obtain $\boldsymbol{v}_d=\mathbf{0}_{t+1}$ or $\mathbf{1}_{t+1}$.
Consequently, $d=1$ or $d=D$ by \eqref{varphi-def}. Thus $\operatorname{Im}(\alpha)=\{1,D\}.$ 

It remains to determine $\operatorname{Im}(\alpha^\prime)$. By Lemma~ \ref{image-torsors} the corresponding quartic equation is
$$C'_d:\quad N^2=dM^4-\frac{4D}{d}e^4.$$ Suppose $C'_d$ has an integer solution $(N,M,e)$ with $\gcd(M,e)=1$. Since $\alpha'((0,0))=-D$ in
$\Q^{*}/\Q^{*2}$, multiplication by $-D$ exchanges positive and negative classes in $\operatorname{Im}(\alpha')$. It therefore suffices to show that its only positive class is $1$. Such a class has representative $d=\delta$ or $2\delta$, where
$\delta\mid D$. Write $D=\delta\delta'$. 

First, suppose that $d=\delta$. The corresponding quartic equation is
\begin{equation}\label{thmB:alpha'}
N^2=\delta M^4-4\delta'e^4.
\end{equation}
In the case considered in Theorem
\ref{TheoremB positive D}, we have
$p_i\equiv1\pmod4$ for all $i$ and $q\equiv3\pmod4$.
Suppose that $q\mid\delta$. Then
$\delta\equiv3\pmod4$ and $\delta'\equiv1\pmod4$.
Reducing \eqref{thmB:alpha'} modulo $4$, we find that $M$ must be even. Since, $\gcd(M,e)=1$ we have $e$
is odd. Reducing the equation modulo $16$, we obtain
 $N^2\equiv-4\delta'\equiv12\pmod{16},$
which is impossible. Therefore, $q\nmid\delta$. Since $\legendre{-1}{p_i}=1$ for every $i$, by Lemma~\ref{rowsumzero for p_i} with $(\lambda, \mu)=(1,-4)$ gives $\mathcal A_D\cdot \boldsymbol{v}_\delta=\mathbf{0}_{t+1}.$
As the last coordinate of $\boldsymbol{v}_\delta$
is zero, the description of $\ker(\mathcal A_D)$
implies that $\boldsymbol{v}_\delta=\mathbf{0}_{t+1}$.
Consequently, $\delta=1$.

Next, suppose that $d=2\delta$. The corresponding quartic equation is
\begin{equation}\label{thmB-alph'-even}
N^2=2\delta M^4-2\delta'e^4.
\end{equation}
Reducing \eqref{thmB-alph'-even} modulo $4$ and using
$\gcd(M,e)=1$, we obtain $N=2N_1$ for some
$N_1\in\mathbb Z$. Now dividing the equation \eqref{thmB-alph'-even} by $2$ we obtain $
2N_1^2=\delta M^4-\delta'e^4.$ Note that $M$ and $e$ are odd, reduction modulo $16$ gives
\[\delta-\delta'\equiv 2N_1^2
\equiv 0,\,2,\text{ or }8\pmod{16}.
\]
If $\delta-\delta'\equiv0$ or $8\pmod{16}$, then
$\delta\delta'\equiv1$ or $9\pmod{16}$.
On the other hand, if $\delta-\delta'\equiv2\pmod{16}$,
then $\delta\delta'\equiv3$ or $15\pmod{16}$.
Each of these possibilities contradicts
$D\equiv7$ or $11\pmod{16}$.
Thus, no positive even class belongs to
$\operatorname{Im}(\alpha')$, and hence $\operatorname{Im}(\alpha')=\{1,-D\}.$ Finally, Proposition~\ref{rank proposition} gives 
$$2^{\operatorname{rank}E_{D}(\Q)+2}
=|\operatorname{Im}(\alpha)| \,|\operatorname{Im}(\alpha')| =4.$$
Consequently, $\operatorname{rank}E_{D}(\Q)=0$. This proves the theorem.
\end{proof}

\begin{proof}[\textbf{\textit{Proof of Theorem \ref{thm-2D-C}}}]
     Let $E=E_{2D}:y^2=x^3+2Dx.$
The elliptic curve $2$-isogenous to $E$ is $E^\prime=E_{-8D}:y^2=x^3-8Dx.$
Let $\alpha$ and $\alpha'$ be the descent homomorphisms defined in \eqref{alpha} and \eqref{alpha'}, respectively. We prove that $\operatorname{Im}(\alpha)=\{1,2D\} \quad \text{and}\quad
\operatorname{Im}(\alpha^\prime)=\{1,-2D\}$ where the entries denote classes in $\Q^*/\Q^{*2}$. All matrix operations below are over $\F_2$. Since every row of $\mathcal A_D$ has sum zero and $\operatorname{rank}_{\F_2}(\mathcal A_D)=t$, we have $\ker(\mathcal A_D)
=\{\mathbf0_{t+1},\mathbf1_{t+1}\}.$ For a positive divisor $\delta\mid D$, write
$\delta=p_1^{e_1}\cdots p_t^{e_t}q^f,
~ e_1,\ldots,e_t,f\in\{0,1\},$
and set $\boldsymbol v_\delta=\varphi(\delta)
=(e_1,\ldots,e_t,f)^{\mathsf T}\in\F_2^{t+1},$
where $\varphi$ is the isomorphism defined in
\eqref{varphi-def}. 

By Lemma~\ref{image-torsors}, every class in either descent image has a square-free representative dividing $2D$.
Moreover,
$$\alpha((0,0))=2D, \quad \alpha'((0,0))=-8D\equiv-2D\pmod{\Q^{*2}}.$$
Both images are subgroups of $\Q^*/\Q^{*2}$ and hence multiplication by the corresponding known class interchanges odd and even classes in each image. It therefore suffices to determine all odd classes (positive or negative).

We first determine $\operatorname{Im}(\alpha)$.
Every class in this image is positive.
Indeed, if a divisor $d$ of $2D$ is negative, then both coefficients in the associated quartic equation
$$N^2=dM^4+\frac{2D}{d}e^4$$ are negative.
Since $\gcd(M,e)=1$ implies that $M,e$ are not both zero, this equation has no such integer solution. Let $\delta\mid D$ be positive and represent an odd class in $\operatorname{Im}(\alpha)$, and put $\delta'=D/\delta$. By Lemma~\ref{image-torsors}, the associated quartic equation
\begin{equation}\label{EQ_3.7}
N^2=\delta M^4+2\delta'e^4
\end{equation}
has an integer solution with $\gcd(M,e)=1$. Suppose that $q\mid\delta$. Since every $p_i\equiv1\pmod8$, we obtain
$\delta\equiv q\pmod8,~ \delta'\equiv1\pmod8.$
If $M$ is even, then $e$ is odd, and \eqref{EQ_3.7} gives $N^2\equiv2\pmod8$.
If $e$ is even, then $M$ is odd, and $N^2\equiv q\pmod8$. If $M,e$ are both odd, then
$N^2\equiv q+2\pmod8$. Since $q\equiv3$ or $5\pmod8$, each possibility contradicts
the fact that the squares modulo $8$ are $0,1,4$. Thus $q\nmid\delta$, and hence $f=0$. Apply Lemma~\ref{rowsumzero for p_i} to \eqref{EQ_3.7} with $(\lambda,\mu)=(1,2)$.
For every $1\le i\le t$, the $i$th coordinate of
$\mathcal A_D\cdot \boldsymbol v_\delta$ is zero:
it is $\left[\frac1{p_i}\right]=0$ if $p_i\nmid\delta$, and $\left[\frac2{p_i}\right]=0$ if $p_i\mid\delta$,
because $p_i\equiv1\pmod8$.
The last coordinate is also zero, since $q\nmid\delta$ and $\left[\frac1q\right]=0$.
Therefore
$$\mathcal A_D\cdot \boldsymbol v_\delta=\mathbf0_{t+1}.$$
Since the last coordinate of $\boldsymbol v_\delta$ is zero, the description of $\ker(\mathcal A_D)$ implies $\boldsymbol v_\delta=\mathbf0_{t+1}$.
Consequently, $\delta=1$.
Thus the only odd class in $\operatorname{Im}(\alpha)$ is $1$. Since $2D$ belongs to this image and multiplication by it interchanges odd and even classes, we conclude that $$\operatorname{Im}(\alpha)=\{1,2D\}.$$

We now determine $\operatorname{Im}(\alpha')$. Let $\delta$ be odd square-free divisor of $D$ representing an odd class in this image, and put $\delta'=D/\delta$.
By Lemma~\ref{image-torsors}, the associated quartic equation
\begin{equation}\label{EQ_3.8}
N^2=\delta M^4-8\delta'e^4
\end{equation}
has an integer solution with $\gcd(M,e)=1$. If $M$ were even, then $e$ would be odd.
Reducing \eqref{EQ_3.8} modulo $16$ would give
$$N^2\equiv-8\delta'\equiv8\pmod{16},$$
which is impossible. Hence $M$ is odd. Equation \eqref{EQ_3.8} then implies that $N$ is odd. Reducing modulo $8$, we obtain $1\equiv N^2\equiv\delta\pmod8.$ Every positive divisor of $D$ is congruent to either $1$ or $q$ modulo $8$. A negative divisor is therefore congruent to either $7$ or $-q$ modulo $8$.
Since $q\equiv3$ or $5\pmod8$, neither negative residue is $1$ modulo $8$. Thus $\delta>0$. Moreover, if $q\mid\delta$, then
$\delta\equiv q\pmod8$, contradicting
$\delta\equiv1\pmod8$. Therefore $q\nmid\delta$, and the last coordinate of
$\boldsymbol v_\delta$ is zero. Apply Lemma~\ref{rowsumzero for p_i} to \eqref{EQ_3.8} with $(\lambda,\mu)=(1,-8)$.
For $1\le i\le t$, we have
$\legendre{-8}{p_i}=1,$
because $p_i\equiv1\pmod8$.
Hence the first $t$ coordinates of
$\mathcal A_D\cdot \boldsymbol v_\delta$ are zero.
The last coordinate is zero because $q\nmid\delta$ and $\left[\frac1q\right]=0$.
Thus $$\mathcal A_D\cdot \boldsymbol v_\delta=\mathbf0_{t+1}. $$
The description of $\ker(\mathcal A_D)$ and the zero last coordinate give $\boldsymbol v_\delta=\mathbf0_{t+1}$.
Therefore $\delta=1$. We have proved that the only odd class in $\operatorname{Im}(\alpha')$ is $1$. Since $-2D$ belongs to this image and multiplication by it interchanges odd and even classes, it follows that $\operatorname{Im}(\alpha')=\{1,-2D\}.$ Finally, Proposition~\ref{rank proposition} gives $\operatorname{rank}E_{2D}(\Q)=0$.
\end{proof}

\begin{proof}[\textbf{\textit{Proof of Theorem \ref{thm-2D-D}}}]
Let $D = p_1 p_2 \cdots p_t$ be a square-free integer, where $p_i \equiv 5 \pmod 8$ for all $i\in \{1, \ldots, t\}$. Consider the elliptic curve over $\Q$ given by $E= E_{2D}: y^2=x^3+2Dx$ and the elliptic curve $2$-isogenous to $E$ is $E'=E_{-8D}:y^2=x^3-8Dx.$
Let $\alpha$ and $\alpha'$ be the corresponding descent maps as defined in \eqref{alpha} and \eqref{alpha'} respectively. It suffices to prove that $$\alpha(E(\Q))=\{1, 2D\} ~ \text{ and } \alpha^\prime(E^\prime(\Q))=\{1, -2D \}.$$ For the curve $E$,  by Lemma~\ref{image-torsors} the associated quartic equation is $$N^2=dM^4+\frac{2D}{d}e^4$$ has an integer solution with $\gcd(M,e)=1$. One can easily see that $d$ must be positive. Since $2D\in\operatorname{Im}(\alpha)$, multiplication
by this class interchanges odd and even classes in $\operatorname{Im}(\alpha)$.
It therefore suffices to determine the odd classes. Let $\displaystyle d=\prod_{i=1}^tp_i^{e_i}$ where each $p_i\equiv 5 \pmod8$, $e_i \in \F_2$ and write $\boldsymbol{v}_{d}=(e_1,\ldots,e_t)^{\mathsf T}$ as in \eqref{varphi-def}. We begin by computing the image of the descent map $\alpha$, namely, $\operatorname{Im}(\alpha)$.

\textbf{Case 1:} $p_i \nmid d$ or equivalently, $e_i=0$. The corresponding quartic equation is $N^2=dM^4+\frac{2D}{d}e^4$. Taking modulo the prime $p_i$ we have $\legendre{d}{p_i}=1$ and consequently, $\displaystyle \sum_{j\neq i}a_{ij}e_j=0$. Note that $$((\mathcal{A}_{D}+\mathcal I_t)\cdot\boldsymbol{v}_{d})_i =\sum_{j=1}^t(a_{ij}+\delta_{ij})e_j=\sum_{j \neq i}a_{ij}e_j+(a_{ii}+1)e_i,$$ where $\delta_{ij}=1$ if  $i=j$, and $\delta_{ij}=0$ if $i \neq j$.  Since $e_i=0$ we have $((\mathcal{A}_{D}+\mathcal I_t)\cdot \boldsymbol{v}_{d})_i =\sum_{j \neq i}a_{ij}e_j=0.$

\textbf{Case 2: } $p_i \mid d$, or equivalently $e_i=1$. Reducing the quartic equation modulo $p_i$ we have $\legendre{2D/d}{p_i}=1$. Since $p_i \equiv 5 \pmod8$, we have $\legendre{2}{p_i}=-1$ and consequently $\left[ \frac{2}{p_i} \right]=1$. Hence we have $1+a_{ii}+\sum_{j\neq i}a_{ij}e_j=0$. Therefore, $((\mathcal{A}_{D}+\mathcal I_t)\cdot \boldsymbol{v}_{d})_i=0$ for all $1\leq i \leq t$.

According to the hypothesis we have that $ \rank_{\F_2}(\mathcal{A}_{D}+\mathcal I_t)=t. $ Hence $\ker(\mathcal{A}_{D}+\mathcal I_t)=\{\mathbf{0}_{t}\}.$ Therefore we have, $\boldsymbol{v}_{d}=\mathbf{0}_{t}$. Hence $d=1.$  Consequently $\operatorname{Im}(\alpha)=\{1,2D\}.$

Next we will calculate $\operatorname{Im}(\alpha^\prime)$.
The isogenous curve is $E^\prime: y^2=x^3-8Dx$. 
 The point $(0,0)\in E'(\Q)$ gives
$\alpha'((0,0))=-8D\equiv-2D\pmod{\Q^{*2}}.$
By Lemma~\ref{image-torsors}, every class in
$\operatorname{Im}(\alpha')$ has a square-free
representative dividing $2D$. Since multiplication by the known class $-2D$ interchanges odd and even classes in this image, it suffices to determine all odd classes, allowing both signs.

Let $d$ be an odd divisor of $D$ representing such a class. By Lemma~\ref{image-torsors}, the associated quartic equation
$$N^2=dM^4-\frac{8D}{d}e^4$$
has an integer solution $(N,M,e)$ with $\gcd(M,e)=1$.
If $M$ were even, then $e$ would be odd, and reduction modulo $16$ would give $N^2\equiv8\pmod{16},$ which is impossible. Hence $M$ is odd, and the equation implies that $N$ is odd. Reducing modulo $8$, we obtain $1\equiv N^2\equiv dM^4\equiv d\pmod8.$
Since every $p_i\equiv5\pmod8$, every positive divisor of $D$ is congruent to either $1$ or $5$ modulo $8$. Every negative divisor is therefore congruent to either $7$ or $3$ modulo $8$. Thus $d\equiv1\pmod8$ forces $d>0$.
Consequently, it remains to determine the positive odd divisors of $D$ representing classes in $\operatorname{Im}(\alpha')$.
For a positive odd divisor $d\mid D$ representing a class
in $\operatorname{Im}(\alpha')$, write
$\displaystyle d=\prod_{j=1}^{t}p_j^{e_j},~ e_j\in\{0,1\},$
and set $\boldsymbol v_d=\varphi(d)
=(e_1,\ldots,e_t)^{\mathsf T}\in\F_2^t,$
as in \eqref{varphi-def}. The associated quartic equation is
$N^2=dM^4-8\frac{D}{d}e^4.$
Since $p_i\equiv5\pmod8$, we have
$\legendre{-8}{p_i}=-1.$
Thus the same local calculation as for $\alpha$ gives
$$(\mathcal A_D\cdot \boldsymbol v_d)_i
=\begin{cases}
0,&p_i\nmid d,\\
1,&p_i\mid d,
\end{cases}
=e_i.$$
Consequently, $\mathcal A_D\cdot \boldsymbol v_d=\boldsymbol v_d,$ i.e. $(\mathcal A_D+\mathcal I_t) \cdot\boldsymbol v_d=\mathbf0_t.$
By the invertibility of the matrix $\mathcal A_D+\mathcal I_t$, we obtain $\boldsymbol v_d=\mathbf0_t$, and hence $d=1$. Therefore the only odd class in
$\operatorname{Im}(\alpha')$ is $1$.
Since $-2D$ belongs to this image and multiplication by it
interchanges odd and even classes, we conclude that $\operatorname{Im}(\alpha')=\{1,-2D\}.$ Finally, Proposition~\ref{rank proposition} gives $\operatorname{rank}E_{2D}(\Q)=0$.
\end{proof}

\begin{proof}[\textbf{\textit{Proof of Theorem \ref{thm 2D-E}}}]
Let $E=E_{2D}:y^2=x^3+2Dx.$ The elliptic curve $2$-isogenous to $E$ is $E'=E_{-8D}:y^2=x^3-8Dx.$
Let $\alpha$ and $\alpha'$ be the corresponding descent maps as defined in \eqref{alpha} and \eqref{alpha'} respectively. We first determine their images.

By Lemma~\ref{image-torsors}, the associated quartic equation $C_d: N^2=dM^4+\frac{2D}{d}e^4$
admits an integer solution $(N,M,e)$ with $\gcd(M,e)=1$ and $d\mid 2D$.
Since $D>0$ and $\alpha((0,0))=2D$, it suffices to consider
positive odd divisors $d$ of $D$.  For a positive $d\mid D$ write $d=p_{1}^{e_{1}}p_{2}^{e_{2}}\cdots p_{t}^{e_{t}},\;\;
\boldsymbol{v}_{d}=(e_1,e_2,\ldots ,e_t)^{\mathsf T}\in\F_2^t$ as in \eqref{varphi-def}. 

If $p_i\nmid d$, then reducing modulo $p_i$ gives $\left(\frac{d}{p_i}\right)=1$, and  hence we have
$\sum_{j\neq i}a_{ij}e_j=0$. Since $e_i=0$, this gives
$(\mathcal{A}_{D}\cdot\boldsymbol{v}_{d})_i=0.$
If $p_i\mid d$, then reducing modulo $p_i$ gives
$\legendre{2D/d}{p_i}=1$. Since $p_i\equiv3\pmod8$, we have
$\legendre{2}{p_i}=-1$, and hence
$$1+\sum_{j=1,j \neq i}^t a_{ij}(1-e_j)=0.$$
Thus $(\mathcal{A}_{D}\cdot\boldsymbol{v}_{d})_i=1$. Combining the two cases, the $i$-th component of $\mathcal{A}_{D}\cdot\boldsymbol{v}_{d}$ equals $0$ when $e_i=0$ and $1$ when $e_i=1$. Hence
$\mathcal{A}_D \cdot \boldsymbol v_d=\boldsymbol v_d,
\quad\text{equivalently}\quad
(\mathcal{A}_D + \mathcal I_t)\cdot\boldsymbol v_d=\mathbf{0}_t.$ Since every $p_i\equiv3\pmod4$, quadratic reciprocity gives
$a_{ij}+a_{ji}=1$ for $i\ne j$.
Moreover, $t$ is odd, the diagonal entries satisfy $\displaystyle a_{ii}=\sum_{j\ne i}a_{ij}$, and
$\operatorname{rank}_{\F_2}(\mathcal A_D)=t-1$. By Lemma \ref{odd lemma invertible} we have $\mathcal A_D+\mathcal I_t$ is invertible, and consequently $\boldsymbol v_d=\mathbf{0}_t$, so $d=1$. Therefore
$\operatorname{Im}(\alpha)=\{1,2D\}.$

We now determine $\operatorname{Im}(\alpha^\prime)$. The point $(0,0)\in E'$ gives the square class $\alpha^\prime((0,0))=-8D \equiv -2D \pmod{\Q^{*2}}.$  Thus it suffices to determine the odd classes in $\operatorname{Im}(\alpha')$.  Let $d$ be an odd square-free divisor of $D$ representing such a class. By Lemma~\ref{image-torsors}, the associated quartic equation $C_d:N^2=dM^4-\frac{8D}{d}e^4$
has an integer solution $(N,M,e)$ with $\gcd(M,e)=1$.

 If $M$ were even, then $\gcd(M,e)=1$ will imply that $e$ is odd. Reducing the quartic equation modulo $16$
would then give $N^2\equiv-8\frac{D}{d}\equiv8\pmod{16},$ which is impossible. Hence $M$ is odd, and the equation
then implies that $N$ is odd. Reducing the quartic equation modulo $8$  gives
$N^2\equiv d\pmod8 \implies d\equiv1\pmod8$. Every positive divisor of $D$ is congruent to either $1$ or $3\pmod{8}$, whereas every negative divisor is congruent to either $5$ or $7\pmod{8}$. Hence, $d$ must be positive. Write
$\displaystyle d=\prod_{j=1}^tp_j^{e_j}, ~ \boldsymbol v_d=(e_1,e_2,\ldots, e_t)^{\mathsf T}$ as in \eqref{varphi-def}.
If $e_i=0$, the local condition
$\left(\frac d{p_i}\right)=1$ gives
$(\mathcal A_D\cdot\boldsymbol v_d)_i=0$.
If $e_i=1$, the local condition is
$\left(\frac{-8D/d}{p_i}\right)=1.$
Because $p_i\equiv3\pmod8$, we have
$\left(\frac{-8}{p_i}\right)=1.$
It follows that
$$0=\sum_{j\ne i}a_{ij}(1-e_j)  =a_{ii}+\sum_{j\ne i}a_{ij}e_j  =(\mathcal A_D \cdot \boldsymbol v_d)_i.$$
Thus $\mathcal A_D\cdot \boldsymbol v_d=\mathbf{0}_t$, and hence
$\boldsymbol v_d=\mathbf{0}_t$ or $\mathbf{1}_t$. Equivalently, $d=1$ or $D$.
But $D\equiv3\pmod8$, whereas $d\equiv1\pmod8$, so $d=1$.
We conclude that
$\operatorname{Im}(\alpha')=\{1,-2D\}.$ Therefore, by Proposition \ref{rank proposition} we are done.
\end{proof}

We have constructed families of elliptic curves $E$ with rank $0$ over $\Q$. In fact, using Proposition \ref{rank over K}, we further show that these families also have rank $0$ over the imaginary quadratic field $\Q(i)$. Let $K = \Q(\sqrt{m})$, where $m$ is a square-free integer. The Mordell-Weil rank of an elliptic curve over $K$ can be expressed in terms of the ranks over $\Q$ of the curve and its $m$-twist. We recall the following result (see, for example, \cite{silverman2009arithmetic}).

\begin{proposition}\label{rank over K}
 Let $K = \Q(\sqrt{m})$ be a quadratic field, where $m$ is a square-free integer. Let $E : y^2 = x^3 + ax^2 + bx$ be an elliptic curve over $\Q$ and $E^m: y^2 = x^3 + max^2 + m^2bx$ be the $m$-twist of $E$. Then 
    $$\operatorname{rank}(E(K)) = \operatorname{rank}(E(\Q)) + \operatorname{rank}(E^m(\Q)).$$    
\end{proposition}

\begin{corollary}\label{allCorollary}
Each family of elliptic curves described in Theorems
\ref{TheoremA negative D}--\ref{thm 2D-E}
has rank $0$ over the imaginary quadratic field $\Q(i)$.
\end{corollary}

\begin{proof}
Recall that the elliptic curves under consideration are of the form $E_c: y^2=x^3+cx$. Since the quadratic twist of $E_c$ by $-1$
is isomorphic to $E_c$ over $\Q$, Proposition \ref{rank over K}
implies that
\[
\operatorname{rank} E_c(\Q(i))
=2\operatorname{rank} E_c(\Q).
\]
By Theorems \ref{TheoremA negative D}--\ref{thm 2D-E},
we have $\operatorname{rank} E_c(\Q)=0$.
Consequently, $\operatorname{rank} E_c(\Q(i))=0$,
as desired.
\end{proof}

\section{Density of rank zero curves}\label{Density section}

In this section, we derive a density consequence of one of our main theorems, namely the rank-zero criterion for the family $E_{2D}:y^2=x^3+2Dx$ in Theorem \ref{thm-2D-C}. We show that a positive lower proportion of the square-free parameters in the families described in that theorem give elliptic curves of Mordell--Weil rank zero over $\Q$. The argument combines the counting of square-free integers with prime factors in prescribed residue classes, the equidistribution of Legendre symbols, and the rank distribution of the corresponding matrices over $\F_2$. Analogous density results can be derived for the remaining main results by adapting a similar approach.
We give the details for this case to illustrate the method and the resulting proportions in the other cases depend on their respective hypotheses.
For real-valued functions $a(x)$ and $b(x)$, we write $a(x)\sim b(x)$ if $\displaystyle \lim_{x\to\infty}\frac{a(x)}{b(x)}=1.$

Let $t \geq 1$ be fixed and let $r\in\{3,5\}$. Define
$$\begin{aligned}
\mathcal D_t^{(r)} &:=\left\{D=p_1\cdots p_tq: p_1,\ldots,p_t,q \text{ are distinct primes},\ p_i\equiv1\pmod8,\ q\equiv r\pmod8 \right\},\\
\mathcal D_t^{(r),\mathrm{rk}}(x) &:=\left\{D\in\mathcal D_t^{(r)}:D\leq x,\  \operatorname{rank}_{\F_2}(\mathcal{A}_{D})=t\right\},\\
\mathcal D_t^{(r),0}(x) &:=\left\{D\in\mathcal D_t^{(r)}:D\leq x,\ 
\operatorname{rank}E_{2D}(\mathbb Q)=0\right\},\\
\pi_k(x)&:=\#\left\{n\leq x:n\text{ is square-free and has exactly }k \text{ distinct prime factors}\right\}.\\
p(t) &:=2^{-t(t+1)/2}\#\left\{B \in M_t(\F_2):
B=B^{\mathsf T},\ \det B \neq 0\right\}.\\
\end{aligned}
$$

\begin{proposition}
With the notation introduced above we have the following:
 $p(t)>0$ and as $x \to \infty$, we have
$$\#\mathcal D_t^{(r),\mathrm{rk}}(x)
\sim\frac{t+1}{4^{t+1}}\,p(t)\,\pi_{t+1}(x), \quad \text{ or equivalently} \quad 
\#\mathcal D_t^{(r),\mathrm{rk}}(x)
\sim \frac{(t+1)p(t)}{4^{t+1}t!}
\frac{x(\log\log x)^t}{\log x}.$$

Consequently, by the rank-zero criterion established in the Theorem \ref{thm-2D-C},
$$\#\mathcal D_t^{(r),0}(x)
\geq \#\mathcal D_t^{(r),\mathrm{rk}}(x)
\sim \frac{t+1}{4^{t+1}}\,p(t)\,\pi_{t+1}(x).$$
\end{proposition}

\begin{proof}
 Put $k:=t+1, \quad m:=\binom{k}{2}=\frac{t(t+1)}{2},$ and write
 $$N_r(x) :=\#\{D \in \mathcal D_t^{(r)}:D\leq x\}.$$

We first count the integers satisfying the prescribed
residue class conditions. For fixed nonnegative integers $a_1,a_3,a_5,a_7$ satisfying $a_1+a_3+a_5+a_7=k,$
let $N(x;a_1,a_3,a_5,a_7)$ denote the number of square-free integers $n\leq x$ with $(n,8)=1$ and exactly $k$ prime factors, of which exactly $a_j$ are congruent to $j \pmod 8$. Rhoades's counting theorem
\cite[Theorem~3.1]{Rhoades}, also cited in
\cite[Proposition~2.7]{Vrecica}, gives
$$N(x;a_1,a_3,a_5,a_7)
\sim
\frac{k!}{a_1!a_3!a_5!a_7!}
\frac{1}{4^k}
\frac{x(\log\log x)^{k-1}}{(k-1)!\log x}.$$

For the family $\mathcal D_t^{(r)}$, we take
$a_1=t$, $a_r=1$, and all remaining multiplicities equal to zero.
Thus
\begin{equation}\label{eq:rank-density-family-count}
N_r(x) \sim \frac{(t+1)!}{t!\,1!} \frac{1}{4^{t+1}}
\frac{x(\log\log x)^t}{t!\log x}
=\frac{t+1}{4^{t+1}}
\frac{x(\log\log x)^t}{t!\log x}.
\end{equation}
Landau's theorem \cite[equation~(8)]{Landau1900} gives
$$\pi_{t+1}(x)\sim \frac{x(\log\log x)^t}{t!\log x}.$$
Consequently, we have
\begin{equation}\label{eq:rank-density-family-ratio}
N_r(x)\sim \frac{t+1}{4^{t+1}}\pi_{t+1}(x).
\end{equation}
For each $D\in\mathcal D_t^{(r)}$, arrange its prime factors in increasing order:
$D=\ell_1\cdots\ell_k, \quad \ell_1<\cdots<\ell_k.$
Changing the labeling of the prime factors simultaneously permutes the rows and columns of $\mathcal A_D$, and therefore preserves its rank. Also note that the matrix $\mathcal{A}_D$ is symmetric. Consider the finite set
$$\mathscr M_t :=\left\{A\in M_k(\F_2):
A=A^{\mathsf T},\ A\cdot \mathbf 1_k=\mathbf{0}_k \right\}.$$
Every $A\in\mathscr M_t$ has a unique representation
\begin{equation}\label{eq:rank-density-block-form}
A=\begin{pmatrix}
B & B\cdot\mathbf 1_t\\
\mathbf 1_t^{\mathsf T}\cdot B
&
\mathbf 1_t^{\mathsf T}\cdot B \cdot\mathbf 1_t
\end{pmatrix},
\quad B=B^{\mathsf T}\in M_t(\mathbb F_2).
\end{equation}
Indeed, the zero row sums determine the last column from $B$, and symmetry determines the last row.
Adding the first $t$ columns of
\eqref{eq:rank-density-block-form} to its last column,
and then adding the first $t$ rows to its last row,
transforms $A$ into
$\begin{pmatrix}
B & 0\\
0 & 0
\end{pmatrix}.
$
All these operations are over $\F_2$ and preserve rank.
Therefore $\operatorname{rank}_{\mathbb F_2}A
=\operatorname{rank}_{\mathbb F_2}B.$
Since a symmetric $t\times t$ matrix has $m=t(t+1)/2$
freely chosen entries, it follows that $\#\mathscr M_t=2^m.$
Writing
$$\mathscr G_t
:=\left\{A\in\mathscr M_t: \operatorname{rank}_{\mathbb F_2}A=t\right\},$$
the definition of $p(t)$ and the preceding bijection give
\begin{equation}\label{eq:rank-density-good-matrices}
\#\mathscr G_t=2^m p(t).
\end{equation}
Taking $B=\mathcal I_t$ shows that $p(t)\geq 2^{-m}>0.$
It remains to establish the limiting distribution of these
matrices among the integers under consideration. For each $s\in\{1,\ldots,k\}$, define
$\delta_i^{(s)}:=r$ if $i=s$, and $\delta_i^{(s)}:=1$ otherwise, 
and put
$$C_s(x):=
\left\{
(\ell_1,\ldots,\ell_k):
\begin{array}{l}
\ell_1,\ldots,\ell_k\text{ are primes},\\
2<\ell_1<\cdots<\ell_k,\\
\ell_1\cdots\ell_k\leq x,\\
\ell_i\equiv\delta_i^{(s)}\pmod 8
\quad(1\leq i\leq k)
\end{array}
\right\}.
$$
Unique prime factorization gives a bijection between the disjoint union of the sets $C_s(x)$ and the integers counted by $N_r(x)$.
Thus
\begin{equation}\label{eq:rank-density-position-sum}
N_r(x)=\sum_{s=1}^{k}\#C_s(x).
\end{equation}

Fix $A=[a_{ij}]\in\mathscr M_t$.
For an increasingly ordered tuple
$(\ell_1,\ldots,\ell_k)\in C_s(x)$,
quadratic reciprocity and the definition of $\mathcal A_D$
give the equivalence
$$\mathcal A_{\ell_1\cdots\ell_k}=A
\quad\Longleftrightarrow\quad
\left(\frac{\ell_i}{\ell_j}\right)=(-1)^{a_{ij}}
\quad\text{for every }i<j.$$
The diagonal entries impose no additional conditions,
because they are determined by the off-diagonal entries and the zero row sums.

Rhoades's Legendre-symbol equidistribution theorem
\cite[Theorem~4.1]{Rhoades}, reproduced in
\cite[Proposition~2.8]{Vrecica},
states that every prescribed pattern of the
$\binom{k}{2}=m$ symbols
$\left(\frac{\ell_i}{\ell_j}\right),
\;\; i<j,$
occurs with limiting proportion $2^{-m}$ among increasingly ordered prime tuples with fixed residue classes modulo $8$ and product at most $x$.
Applying this theorem to $C_s(x)$ yields
\begin{equation}\label{eq:rank-density-matrix-distribution}
\#\left\{
(\ell_1,\ldots,\ell_k)\in C_s(x):
\mathcal A_{\ell_1\cdots\ell_k}=A
\right\}
=2^{-m}\#C_s(x)+o\bigl(\#C_s(x)\bigr).
\end{equation}
Here $o(\#C_s(x))$ denotes an error term that is negligible compared with $\#C_s(x)$ as $x\to\infty$. Summing \eqref{eq:rank-density-matrix-distribution}
over $A\in\mathscr G_t$ and $1\leq s\leq k$, we obtain
$$\begin{aligned}
\#\mathcal D_t^{(r),\mathrm{rk}}(x)
&=\sum_{s=1}^{k}\sum_{A\in\mathscr G_t}
\left(
2^{-m}\#C_s(x)+ o\bigl(\#C_s(x)\bigr)
\right)\\
&=2^{-m}\#\mathscr G_t
\sum_{s=1}^{k}\#C_s(x)+ o\bigl(N_r(x)\bigr)=p(t)N_r(x)+o\bigl(N_r(x)\bigr),
\end{aligned}
$$
where we used \eqref{eq:rank-density-good-matrices} and \eqref{eq:rank-density-position-sum}.
The summation of the error terms is valid because $k$ is fixed and $\mathscr G_t$ is a fixed finite set. Since $p(t)>0$, combining the last equation with
\eqref{eq:rank-density-family-ratio} gives
$$\#\mathcal D_t^{(r),\mathrm{rk}}(x)
\sim \frac{t+1}{4^{t+1}}\,p(t)\,\pi_{t+1}(x).$$

Combining it instead with
\eqref{eq:rank-density-family-count} gives the equivalent asymptotic
$\#\mathcal D_t^{(r),\mathrm{rk}}(x)
\sim
\frac{(t+1)p(t)}{4^{t+1}t!}
\frac{x(\log\log x)^t}{\log x}.$ Finally, the rank-zero criterion from Theorem~\ref{thm-2D-C} gives
$\operatorname{rank}_{\mathbb F_2}\mathcal A_D=t
\Longrightarrow\quad
\operatorname{rank}E_{2D}(\mathbb Q)=0.$
Therefore
$$\mathcal D_t^{(r),\mathrm{rk}}(x)
\subseteq \mathcal D_t^{(r),0}(x),$$
which proves the asserted lower bound.

\end{proof}

\section{Numerical Computation}\label{Numerical Section}

In this section, we present numerical examples illustrating Theorems~ \ref{TheoremA negative D}--\ref{thm 2D-E}. For each theorem, the examples are grouped according to its cases, with the chosen primes satisfying the stated hypotheses. The tables summarize the relevant primes, the parameter $D$, the rank over $\F_2$ of the associated matrix, and the Mordell--Weil rank of the corresponding elliptic curve over $\Q$. All computations were performed using SageMath~\cite{sage}. 

\begin{table}[htbp]
\centering
\label{tab:minusD}
\small
\begin{tabular}{l@{\hspace{0.3cm}}l@{\hspace{0.7cm}}l@{\hspace{0.7cm}}l@{\hspace{0.7cm}}l@{\hspace{0.7cm}}l@{\hspace{0.7cm}}l}
\toprule
Case & $t$ & $(p_1,\ldots,p_t)$ & $q$ & $D$
& $\rank_{\F_2}(\mathcal{A}_{D})$
& $\rank E_{-D}(\Q)$ \\
\midrule
(i)  & $1$ & $(17)$ & $3$ & $17\cdot 3$ & $1$ & $0$ \\
(i)  & $2$ & $(17,\ 41)$ & $3$ & $17\cdot 41\cdot 3$ & $2$ & $0$ \\
(i)  & $3$ & $(17,\ 41,\ 73)$ & $3$ & $17\cdot 41\cdot 73\cdot 3$ & $3$ & $0$ \\
\midrule
(ii) & $1$ & $(5)$ & $7$ & $5\cdot 7$ & $1$ & $0$ \\
(ii) & $1$ & $(13)$ & $71$ & $13\cdot 71$ & $1$ & $0$ \\
(ii) & $3$ & $(5,\ 13, \ 37)$ & $71$ & $5\cdot 13 \cdot 37 \cdot 71$& $3$ & $0$ \\
\midrule
(iii) & $2$ & $(5,\ 13)$ & $3$ & $5\cdot 13\cdot 3$ & $2$ & $0$ \\
(iii) & $2$ & $(5,\ 13)$ & $11$ & $5\cdot 13\cdot 11$ & $2$ & $0$ \\
\midrule
(iv) & $1$ & $(17)$ & $29$ & $17 \cdot 29$ & $1$ & $0$ \\
(iv) & $2$ & $(17,\ 41)$ & $5$ & $17\cdot 41 \cdot 5$ & $2$ & $0$\\
(iv) & $3$ & $(17, \ 41 ,\ 73)$ & $61$ & $17\cdot 41 \cdot 73 \cdot 61$ & $3$ & $0$ \\
\bottomrule
\end{tabular}
\par\vspace{5pt}
\caption{Numerical examples for Theorem~\ref{TheoremA negative D}.}
\end{table}

\begin{table}[htbp]
\centering
\small
\setlength{\tabcolsep}{5pt}
\begin{tabular}{@{}cclrlcc@{}}
\toprule
$D \bmod 16$ & $t$ & $(p_1,\ldots,p_t)$ & $q$ & $D$
& $\rank_{\F_2}(\mathcal{A}_{D})$
& $\rank E_D(\Q)$ \\
\midrule
$7$ & $1$ & $(17)$  & $7$  & $17\cdot 7$  & $1$ & $0$ \\
$7$ & $2$ & $(17,41)$ & $31$ & $17\cdot 41\cdot 31$  & $2$ & $0$ \\

$7$ & $3$ & $(5,13,29)$  & $3$  & $5\cdot 13\cdot 29\cdot 3$  & $3$ & $0$ \\
\midrule
$11$ & $1$ & $(17)$  & $107$ & $17\cdot 107$ & $1$ & $0$ \\
$11$ & $2$ & $(5,13)$  & $11$ & $5\cdot 13\cdot 11$ & $2$ & $0$ \\
$11$ & $3$ & $(17,41,73)$  & $43$   & $17\cdot 41\cdot 73\cdot 43$  & $3$ & $0$ \\
\bottomrule
\end{tabular}
\par\vspace{5pt}
\caption{Numerical examples for Theorem~\ref{TheoremB positive D}.}
\label{tab:plusD}
\end{table}

\begin{table}[htbp]
\centering
\label{tab:2D-one-five}
\small
\begin{tabular}{l@{\hspace{0.3cm}}l@{\hspace{0.7cm}}l@{\hspace{0.7cm}}l@{\hspace{0.7cm}}l@{\hspace{0.7cm}}l@{\hspace{0.7cm}}l}
\toprule
$t$ & $(p_1,\ldots,p_t)$ & $q$ & $D$
& $\rank_{\F_2}(\mathcal{A}_{D})$
& $\rank E_{2D}(\Q)$ \\
\midrule
 $2$ & $(17, \ 41)$ & $3$ & $17 \cdot 41 \cdot 3$ & $2$ & $0$ \\
 $2$ & $(17, \ 41)$ & $11$ & $17 \cdot 41 \cdot 11$ & $2$ & $0$ \\
 $2$ & $(17,\ 41)$ & $5$ & $17\cdot 41 \cdot 5$& $2$ & $0$ \\
 $3$ & $(17,\ 41, \ 73)$ & $5$ & $17 \cdot 41 \cdot 73 \cdot 5$ & $3$  & $0$ \\
 $3$ & $(17,\ 41, \ 73)$ & $37$ & $17 \cdot 41 \cdot 73 \cdot 37$ & $3$  & $0$ \\
 $3$ & $(17,\ 41, \ 89)$ & $5$ & $17 \cdot 41 \cdot 89 \cdot 5$ & $3$  & $0$ \\
\bottomrule
\end{tabular}
\par\vspace{5pt}
\caption{Numerical examples for Theorem~\ref{thm-2D-C}.}
\end{table}

\begin{table}[htbp]
\centering
\label{tab:2D-five}
\small
\begin{tabular}{l@{\hspace{0.3cm}}l@{\hspace{0.7cm}}l@{\hspace{0.7cm}}l@{\hspace{0.7cm}}l@{\hspace{0.7cm}}l@{\hspace{0.7cm}}l}
\toprule
$t$ & $(p_1,\ldots,p_t)$ & $D$
& $\rank_{\F_2}(\mathcal{A}_{D})$
& $\rank_{\F_2}(\mathcal{A}_{D}+\mathcal I_t)$
& $\rank E_{2D}(\Q)$\\
\midrule
3 &$(5,\ 13,\ 29)$ & $5 \cdot 13 \cdot 29$ & $1$ & $3$ &$0$ \\
3 &$(5,\ 13,\ 61)$ & $5 \cdot 13 \cdot 61$ & $1$ & $3$ &$0$ \\
3 &$(5,\ 29,\ 53)$ & $5 \cdot 29 \cdot 53$ & $1$ & $3$ &$0$ \\
\bottomrule
\end{tabular}
\par\vspace{5pt}
\caption{Numerical examples for Theorem~\ref{thm-2D-D}.}
\end{table}

\begin{table}[htbp]
\centering
\label{tab:three-primes}
\small
\begin{tabular}{l@{\hspace{0.3cm}}l@{\hspace{0.7cm}}l@{\hspace{0.7cm}}l@{\hspace{0.7cm}}l@{\hspace{0.7cm}}l@{\hspace{0.7cm}}l}
\toprule
$p_1$ & $p_2$ & $p_3$ & $D$
& $\rank_{\F_2}(\mathcal{A}_{D})$
& $\rank E_{2D}(\Q)$\\
\midrule
$3$ & $11$ & $19$ & $3\cdot 11 \cdot 19$ & $2$  & $0$ \\
$3$ & $43$ & $59$ & $3\cdot 43 \cdot 59$ & $2$  & $0$ \\
$11$ & $19$ & $67$ & $11\cdot 19 \cdot 67$ & $2$  & $0$ \\
\bottomrule
\end{tabular}
\par\vspace{5pt}
\caption{Numerical examples for Theorem~\ref{thm 2D-E}.}
\end{table}

\newpage
\subsection*{\it Acknowledgement} The first author acknowledges the support of the DST-INSPIRE Faculty Fellowship [DST/INSPIRE/04/2024/004189]. The second author acknowledges the Institute Fellowship and the excellent research environment provided by the Indian Institute of Technology Madras. The third author acknowledges the support of the CSIR-UGC senior research fellowship.\\

\noindent \textbf{Declarations}: \textbf{Conflicts of Interest.} The authors declare that they have no conflict of interest or any competing financial interests that could influence the work reported in this paper.

\bibliographystyle{alpha}
\bibliography{references}

\end{document}